\documentclass[12pt,amsart]{iopart}
\usepackage{fullpage,graphicx,epsf}
\usepackage{iopams}
\usepackage{amssymb}
\usepackage{amscd}
\usepackage{color}
\def\prob{{\rm P}}
\def\eN{{N}}
\def\kK{{R}}
\def\mM{{Q}}
\begin{document}
\title{Matrix representation of the stationary measure\\ for the
multispecies TASEP} 
\author{Martin R.  Evans$^{1,2}$, Pablo A. Ferrari$^{3}$ and Kirone
  Mallick$^{4}$}
\address{
$^1$ SUPA, School of Physics, University of Edinburgh, \\
Mayfield Road, Edinburgh EH9 3JZ, UK.\\
$^2$ Laboratoire de Physique Th\'eorique et Mod\`eles
Statistiques,\\ Universit\'e Paris-Sud, Bat 100, 91405, Orsay Cedex,
France.\\
$^3$ Instituto de Matem\'atica e Estat\'{\i}stica,  Universidade de S\~ao Paulo\\
Caixa Postal 66281
05311-970 S\~ao Paulo,
Brasil\\
$^4$ Institut de Physique Th\'eorique , CEA Saclay, \\ 91191, Gif-sur-Yvette
Cedex, France.
\ead{m.evans@ed.ac.uk,pablo@ime.usp.br,kirone.mallick@cea.fr}
}

\paragraph{\bf Abstract}
In this work we construct the stationary measure of the $N$ species totally
asymmetric simple exclusion process in a matrix product formulation.  We make
the connection between the matrix product formulation and the queueing theory
picture of Ferrari and Martin. In particular, in the standard representation,
the matrices act on the space of queue lengths. For $N >2$ the matrices in
fact become tensor products of elements of quadratic algebras. This enables us
to give a purely algebraic proof of the stationary measure which we present for
$N =3$.  \medskip



\section{Introduction}
Models of diffusing particles with hard core interactions were first
considered in the mathematical literature \cite{Harris65} and the name
exclusion process was first coined by Spitzer \cite{Spitzer70}. In the
totally asymmetric simple exclusion process (TASEP) particles jump
only to the right on a one-dimensional lattice but cannot occupy the
same site. Mathematical achievements include categorising the
stationary measures for the process on ${\mathbb Z}$ and many results
are summarised in the books by Liggett \cite{Liggett85,Liggett99}.

Since the early 1990s the TASEP has been of considerable interest within the
physics community as a prototypical model of nonequilibrium behaviour where,
in the steady state, a current of particles is supported. In particular the
model has been studied on the ring ${\mathbb Z}_L$ and also on a lattice of
length $L$ with open boundary conditions where particles enter at the left
boundary and leave at the right boundary. Notable achievements have been the
use of the Bethe ansatz to determine spectral properties of the transition
rate matrix \cite{Dhar,GS,GM06,dGE06} and the determination of the stationary
state in the open boundary case within a matrix product formulation
\cite{DEHP93}.

A generalisation of the TASEP is to the case of several species of particle.  In
the two-species exclusion process containing first-class particles and
second-class particles \cite{DJLS93} both first and second-class particles hop
to the right with rate 1. However if the site to the right of a first-class
particle is occupied by a second-class particle the first and second-class
particle exchange places with rate 1. Thus a second-class particle behaves as a
hole from the point of view of a first-class particle but behaves as a particle
from the point of view of a hole. The introduction of such a second class
particle is a useful tool to study the microscopic structure of shocks
\cite{FKS,Pablo91,DJLS93}.  Besides, the second-class particle problem arises
naturally from coupling two TASEPs with different densities of particles
\cite{Liggett76}: the excess particles in the system with higher density acquire
the dynamics of second-class particles. The stationary state of a system
containing second and first-class particles has been obtained using the matrix
product formulation by Derrida {\it et al.} \cite{DJLS93}.  Based on this work,
Ferrari, Fontes and Kohayakawa \cite{FFK94} introduced a probabilistic
construction of the measure. Angel \cite{Angel06} improved this construction
providing a combinatorial description of the stationary state.  In \cite{pablo2,FM07}, Ferrari and Martin showed that Angel's work could be interpreted as a
queueing system and they generalized it to the $\eN$ species case, for arbitrary
$\eN$.

  In the physics literature, the exclusion process with $\eN$ species of
  particles was  considered by Mallick, Mallick and Rajewsky
  \cite{MMR99} and studied for the case $\eN =3$.  This model, which we refer to
  as the $\eN$-TASEP, is defined by having site variables $\tau_i$ which may
  take values $0,1\ldots, N$ where $\eN$ is the number of species. (Note that one
  could alternatively consider the state $\tau_i=0$ (a hole) as a species which
  would imply a total of $\eN +1$ species; we choose instead to use the more
  common convention.)  The dynamics is defined as follows: each bond between
  neighbouring lattice sites has a bell which rings with rate 1. When the bell
  at bond $i, i+1$ rings the site variables at $i$ and $i+1$ are exchanged
  provided $\tau_{i+1}=0, \tau_i >0$ or $\tau_{i+1}> \tau_{i} \geq 1$.  This is
  equivalent to the following exchanges occurring with rate 1
\begin{eqnarray}
K \ 0 \to 0\ K&&\qquad \mbox{for}\quad \eN  \geq K\geq 1
\label{eq:rules1}\\
K\ J  \to J \ K &&\qquad \mbox{for}\quad  \eN  \geq J > K \geq 1\;.
\label{eq:rules2}
\end{eqnarray}
The construction of Ferrari and Martin \cite{pablo2,FM07} couples $\eN$
realizations of the TASEP in a special way, called the $\eN$-line process.  To
the $\eN$ configurations in the $\eN$-line system one associates a configuration
of the $\eN$-TASEP. Furthermore, each dynamical event of the $\eN$-line process
corresponds precisely to a dynamical event in the $\eN$-TASEP.  The steady state
measure of the $\eN$-line system is just a uniform distribution of particles.
This implies that 
  one may sample  the $\eN$-TASEP configurations  with their 
stationary state probability  by a two step procedure: (a) uniformly sampling a
configuration of the $\eN$-line system of particles and (b) finding the
associated configuration of the $\eN$-TASEP.

Our aim in this work is to invert this construction to obtain direct expressions
for the steady state probabilities which generalise those already obtained for
the two species case \cite{DJLS93} and the three species case \cite{MMR99}. In
doing so we shall see how the matrix product formulation generalises into a
tensor product.

The paper is structured as follows. In section 2 we review the known results on
the stationary measure of the 2-TASEP and show how the matrix product
representation   \cite{DJLS93} is related to
the queueing representation 
\cite{FM07}. In Section 3 we consider the $\eN$-TASEP and deduce a procedure for
computing the stationary state probabilities. In section 4 we construct a matrix
product representation of the 3-TASEP stationary measure. In section 5 we show
how matrix product representations of the $\eN$-TASEP may be obtained
recursively and we conclude in section 6.

\section{Two Species TASEP}
\label{sec:2}
In this section, we review the known solution of the two species TASEP. We also
illustrate the equivalence between the matrix product solution of Derrida et
al., the construction of Angel and the queueing process interpretation of
Ferrari and Martin.
\begin{figure}[th]
\begin{center}
  \includegraphics[height=5.0cm]{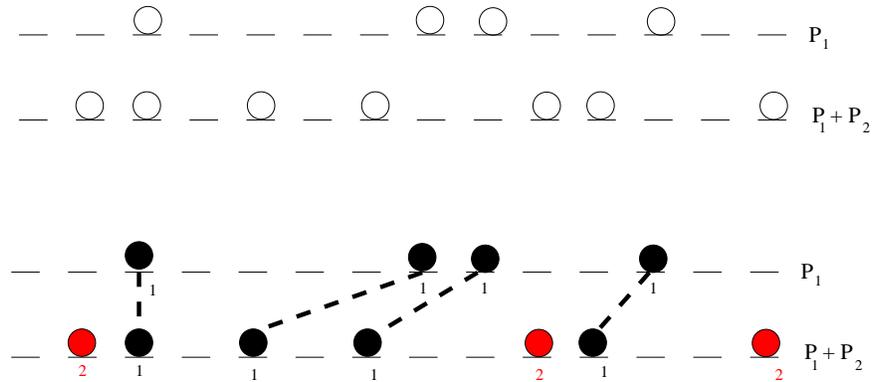}
   \caption{Graphical representation of the construction
for two species.}
  \label{fig:2s}
\end{center}
\end{figure}

\subsection{Construction of Angel that generates the Stationary State}
\label{sec:angel}
We begin by considering the construction of Angel for the two species TASEP on
the ring ${\mathbb Z}_L$. (Note that we often use a different notation to
\cite{Angel06} in order to avoid a clash with some standard notation from the
matrix product formulation.) The construction is to consider a two-line
configuration of particles (see Figure \ref{fig:2s}).  On line 1 there are $P_1$
particles distributed randomly (with at most one particle per site) and on line
2 there are $P_1+P_2$ particles distributed randomly.  Working from right to
left we associate to each particle in line 1, the nearest particle, at the same
site or to the left, in line 2 that has not been associated to another particle.
The associated particles in line 2 are then labelled 1 and the remaining $P_2$
unassociated particles are labelled 2. The empty sites of line 2 are labelled 0
thus each site of line 2 is labelled 0, 1 or 2.  In this way a configuration of
the two species TASEP containing $P_1$ first-class and $P_2$ second-class
particles has been generated through the construction.  Since we consider
periodic boundary conditions, the site at which we begin this procedure (chosen
as the furthest site to the right in Figure \ref{fig:2s}) does not affect the
two species TASEP configuration that is generated, but the particular particle
in line 2 associated to a given particle in line 1 may depend on the initial
starting point; for instance, if in Figure \ref{fig:2s} one starts with the
third particle from the right in line 1, then it would be associated to the
fourth particle from the right in line 2, while the second particle from the
right in line 1 would be associated to the fifth particle from the right in line
2.  As noted in the introduction, uniformly sampling the 2-line configurations
generates 2-TASEP configurations according to their stationary measure.

Angel showed that by uniformly sampling the
two-line configurations, the configurations of the two-species TASEP,
which we denote ${\cal C}$, are sampled with the following
probabilities.
\begin{equation}
\prob({\cal C}) = \frac{\prod_{j=1}^{P_2} \omega(B_j)}{Z(L,P_1,P_2)}
\label{P2}
\end{equation}
Here $\omega(B_j)$ is the weight of the binary string $B_j$ of 0,1
 separating the second-class particles indexed by $j$ and $j+1$
  (here $j\in{\mathbb Z}_{P_2}$ indexes the second class particles). The
normalization
\begin{equation}
\label{Z1}
Z(L,P_1,P_2) =
{{L}\choose{P_1}}{{L}\choose{P_1+P_2}}
\end{equation}
just counts the number of possible 2-line configurations.  Note that the form of
the measure (\ref{P2}) implies a factorization of the stationary state about the
positions of the second-class particles.  The reason for the factorisation is,
as can be seen from Figure \ref{fig:2s}, that all 2-line configurations
associated with a given 2-TASEP configuration must have the following
properties: consider a site $i$, such that there is a particle labelled 2 at $i$
in line 2, then $i$ must be empty on line 1; moreover, no particle in line 1 to
the right of $i$ can be associated to a particle in line 2 to the left of
$i$. This factorization property, which appeared in the matrix product
formulation of Derrida et al \cite{DJLS93}, was used in the construction of the
stationary weights by \cite{FFK94}.

 The weights $\omega(B)$ are given
by the following algorithm which we shall refer to as the {\it pushing procedure:}
given the binary string $B$, one enumerates the number of strings which can be
obtained from it by pushing the 1s to the right, in addition to the original
string. For example from the string $110$ one obtains $110$, $101$, $011$. Thus
$\omega(110)=3$. Similarly, one can obtain from $1010$ the strings $1010$, $0110$,
$1001$, $0101$, $0011$. Thus $\omega(1010)=5$.

 The measure given by (\ref{P2}) is stationary under the dynamics
  of the 2 species TASEP \cite{DJLS93,FFK94,Angel06,Speer94}.  Two key
  properties of this measure are i) the factorisation of the probabilities of
  the 2-TASEP configurations about the position of the second class particles ii)
  the weights $\omega(B)$ are given by the pushing procedure described above.  

  \subsection{Matrix Product Solution of Derrida et al.}
\label{sec:mp2}
 The matrix product formulation has been used to write down the stationary
probabilities of various interacting particle models, thus allowing 
models to be solved through the explicit computation of physical quantities
of interest such as currents, density profiles, correlation functions.
It was first used to solve the TASEP on a lattice of length $L$ with open
boundary  conditions \cite{DEHP93}. It has been extended to 
the 2-TASEP on the ring ${\mathbb Z}_L$ \cite{DJLS93}, partially asymmetric
processes and more general reaction-diffusion systems (for a review see \cite{BE07}).
In this matrix product formulation properties of the stationary measure
manifest themselves in algebraic relations amongst the matrices involved.
Some of these relations have been classified as quadratic algebras 
\cite{IPR01}. 

It is important to note that the measure (\ref{P2}), along with the calculation
of the weights $\omega(B)$, is equivalent to that first obtained within the
matrix product approach, as we now show.  We recall that we use the variable
$\tau_i =0,1,2$ which implies that site $i$ is empty, contains a first-class
particle or contains a second-class particle, respectively.
  Let us denote by 
  ${\cal C}=(\tau_1,\dots,\tau_L)$, a  configuration of the system. In the matrix
product formulation \cite{DJLS93} it has been proved that the stationary measure
may be written as
\begin{equation}
  \prob( {\cal C}) =
  Z^{-1} W({\cal C})\;,
\label{P2mat}
\end{equation}
where the weight  of the configuration is given by
\begin{equation}
W({\cal C})= \mbox{Tr} \left[
    \prod_{i=1}^L X_{\tau_i} \right]
\end{equation}
and Tr means the trace of the product of matrices $X_{\tau_i}$.  The
normalization  $Z=Z(L,P_1,P_2)$  is chosen so that the sum of all the
probabilities is equal to 1.  The matrices $X_{\tau_i}$ are given by
\begin{equation}
X_0= E\;,\quad X_1=D\;,\quad X_2= A\;,
\end{equation}
that is:
if the site is empty we write a matrix $E$;
if the site contains a first-class particle  we write a matrix $D$;
if the site contains a second-class particle  we write a matrix $A$.

The matrices $D,E,A$ obey the algebraic rules
\begin{eqnarray}
DE &=& D + E \label{DE}\\
DA &=& A\label{DA}\\
AE &=& A \; . \label{AE}
\end{eqnarray}
The only remaining condition to satisfy is that representations of $E$,$D$,$A$
may be found which give well-defined values for the traces appearing in
(\ref{P2mat}).  This may be achieved as follows. Let $|n\rangle$ and $\langle
n|$ be the column, respectively row, vector having a 1 in the $n$th coordinate
and 0 in the other ones, $n=0,1,2,\dots\;$. 
Let $A$ be the projector matrix
\begin{equation}
A = |0\rangle \langle 0|
\label{A}
\end{equation}
then $D$,$E$ may be chosen to be bidiagonal semi-infinite matrices
\begin{eqnarray}
D &=& \sum_{n=0}^{\infty}  |n\rangle \langle n| + |n\rangle \langle n+1|
\label{D}\\
E &=& \sum_{n=0}^{\infty}  |n\rangle \langle n| + |n+1\rangle \langle n| \;.
\label{E}
\end{eqnarray}
Writing out the matrices explicitly we have
\begin{eqnarray}
D=               \left( \begin{array}{ccccc}    
                  1&1&0&0&\dots\\
                  0&1&1&0&\ddots\\
                  0&0&1&1&\ddots\\
                  0&0&0&1&\ddots\\              
                  \vdots &\ddots&\ddots&\ddots&\ddots
                        \end{array}  \right)\;,
\qquad
E=
              \left( \begin{array}{ccccc}    
                  1&0&0&0&\dots\\
                  1&1&0&0&\ddots\\
                  0&1&1&0&\ddots\\
                  0&0&1&1&\ddots\\
                  \vdots &\ddots&\ddots&\ddots&\ddots
                        \end{array}  \right) \;,
\\[1ex]
A=
              \left( \begin{array}{ccccc}    
                  1&0&0&0&\dots\\
                  0&0&0&0&\ddots\\
                  0&0&0&0&\ddots\\
                  0&0&0&0&\ddots\\
                  \vdots &\ddots&\ddots&\ddots&\ddots
                        \end{array}  \right)\;.
\end{eqnarray}


Due to the form of $A$, (\ref {P2mat}) reduces to
\begin{equation}
\prob(\left\{ \tau_i \right\}) = Z^{-1} \prod_{j=1}^{P_2}
\omega(B_j)
\label{P2mat2}
\end{equation}
where $B_j$ is as in (\ref{P2}), and $\omega(B_j)$ is now given by
\begin{equation}
\omega(B_j) = \langle0 | \prod_{i = 1}^l X_{\tau_i} | 0 \rangle
\end{equation}
where $l$ is the length of the binary string $B_j$ and $i$ labels the entries in
that string; $X_{\tau_i}$ is either a matrix $D$ or a matrix $E$ according to
whether the entry $\tau_i$ in the string is 1 or 0.

The algebraic rules (\ref{DA},\ref{AE}) imply immediately
that the weight of a string comprising a segment of consecutive zeros followed
by a segment of consecutive ones is equal to 1. In other words, a string where
any 0s are all to the left and any  1s are all to the right has weight 1:
\begin{equation}
\omega(0\cdots 01 \cdots1) = 
\langle0 |E\cdots ED\cdots D | 0 \rangle =
\langle0 | 0 \rangle = 1\;.
\end{equation}
Using rule (\ref{DE}) all binary strings can be reduced to strings of the 
above type and the weight of any string is easily computed. For example,
\begin{eqnarray}
&\omega(10) = \langle0 | DE| 0 \rangle =
\langle0 | D| 0 \rangle 
+\langle0 | E| 0 \rangle   =2\\
&\omega(110) = 
\langle0 |D DE| 0 \rangle =
\langle0 |D D| 0 \rangle 
+\langle0 |D E| 0 \rangle=
1 + 2 =3\nonumber \\
&\omega(1010) = 
\langle0 |DE DE| 0 \rangle =
\langle0 |DDE| 0 \rangle 
+\langle0 |ED E| 0 \rangle
 =
3 + \langle0 |DE | 0 \rangle  =5\nonumber\;.
\end{eqnarray}
  This reduction procedure gives precisely
the same result as the pushing procedure of Angel. In fact Lemma 2.3 of
\cite{FFK94} proves that when $\omega$ is defined via the pushing procedure, the
following relation holds
  \begin{equation}
    \label{p1}
    \omega(B10B') = \omega(B1B')+\omega(B0B')
  \end{equation}
  for arbitrary finite binary sequences $B,\,B'$. But this is the same reduction
  formula that holds for the matrix representation:
\begin{eqnarray}
  \label{p1prime}
  \omega(B10B') &=& \langle0 |XDE X'| 0 \rangle \\
&=&\langle0 |XD X'|
0 \rangle+ \langle0 |XE X'| 0 \rangle
\end{eqnarray}
where $X$, $X'$ are the matrix representation of the binary sequences $B$, $B'$,
respectively. This shows that the definition of $\omega$ 
 by the  matrix  formulation and that by  the pushing procedure coincide.

 The weights $W$ of 2-TASEP configurations are computed from the weights
 $\omega$ of binary strings as follows. Recalling that $N$-TASEP configurations
 are translationally invariant under the periodic boundary conditions we have:
\begin{eqnarray}
&W(0210)\, =\, W(2100)\, =\, \omega(100)\, =\, 3
\label{0210},\\
&W(0211021)\, =\, W(1021102)\, =\, \omega(10)\omega(110) \, =\, 6,\label{6Z}
\end{eqnarray}
and the corresponding probabilities are given by
\begin{eqnarray}
&\prob(0210)\, =\, \prob(2100)\, =\, \frac{W(0210)}{Z(4,1,1)}\, =\, \frac3{24}
\label{0210a},\\
&\prob(0211021)\, =\, \prob(1021102)\, =\, \frac{W(1021102)}{Z(7,3,2)} \, =\, \frac6{735},\label{6Za}
\end{eqnarray}
where $Z$ is defined in (\ref{Z1}).

\subsection{Queueing Interpretation of Ferrari and Martin}
\label{2sq} 
Ferrari and Martin used $\eN$-line configurations to generate $\eN$-TASEP
configurations in terms of queueing processes. Here, we recall the
interpretation in the case of the 2-species TASEP in terms of queueing processes
and make the connection with the matrix product representation of the stationary
measure.

 \begin{figure}[th]
\begin{center}
  \includegraphics[height=5.0cm]{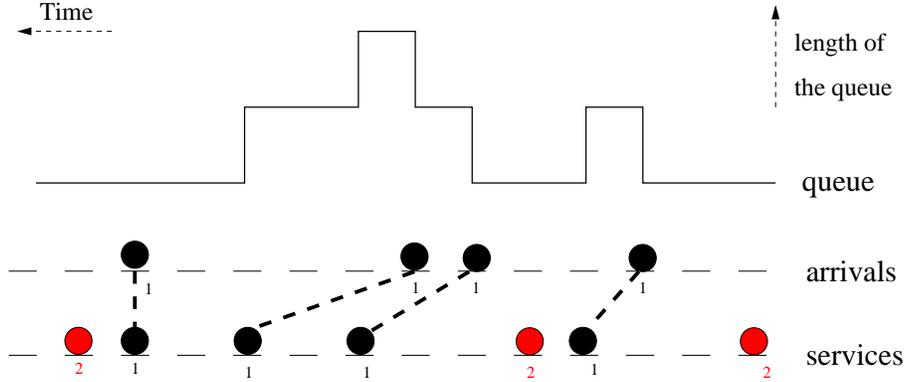}
  \caption{Queueing picture for the two species case. Here time
      $t(i)=L-i$ increases from right to left, where $i$ is site number.}
  \label{fig:2sq}
\end{center}
\end{figure}

We recall (see e.g.~Figure~\ref{fig:2sq}) that a 2-line configuration generates
a 2-TASEP configuration. Since, as described above, the stationary state
factorises about the positions of the second-class particles we only need to
consider a binary string $B$ of 1s and 0s between two 2s in a 2-TASEP
configuration.  There are several 2-line configurations which generate the
string $B$.  Those 2-line configurations must satisfy two conditions: (a) both
lines 1 and 2 must contain the same number of particles, equal to the number of
1s in $B$ and (b) line 2 coincides with $B$.  The possibilities for line 1 are
then generated from line 2 by pushing particles to the right. For example, the
lower line of Figure~\ref{fig:2sq} has three strings of type $B$ delimited by
the three second class particles: $B_0=0$, $B_1=1010100$ and $B_2
  =100$. The upper-line string $1000011$ is one of the strings producing $B_1$,
  the string $010$ is one of the strings  producing $B_2$  and the string $0$ is the
  only string producing $B_0$. 

     Given a 2-line configuration one can associate to it the
    trajectory of the length of a queue. Consider the labels of the particles of
    line 2 as first or second class particles as given by Angel's algorithm
    (illustrated in Figure \ref{fig:2s}). Time for the queue runs from right to
    left: at each site $i$ it is assigned a time $t(i)=L-i$. The queue has length
    zero at the times corresponding to the positions of second class particles
    in line 2 (unused service times), a particle in line 1 represents an arrival
    time and a first class particle in line 2 represents a service time. At a
    given time $t(i)$ the length of the queue (constrained to be non-negative)
    increases by one when a particle is present at site $i$ in line 1 but not in
    line 2 (a new arrival occurs and is not serviced); the length of the queue
  decreases by one when a particle is present in line 2 but not in line 1 (a
  service occurs with no new arrival). If no particles are present in lines 1
  and 2 (no service or new arrival occurs) or when particles are present in both
  lines 1 and 2 (a new arrival occurs and is serviced) the queue remains at the
  same length.  The weight of a 2-TASEP string is then given by all possible
  queue trajectories, consistent with following constraints {\it i}) the queue
  has length zero at the positions of the second class particles {\it ii}) the
  effective service times of the queue are fixed by the positions of the first
  class particles. Since the full 2-line configuration can be retrieved by
  knowing the 2-TASEP configuration and the trajectory of the queue, to
  enumerate the  ancestors of a 2-TASEP configuration it is enough
  to  enumerate the  queue trajectories compatible with it.

We now illustrate how the product of matrices $A$, $D$, $E$ defined in
(\ref{A},\ref{D},\ref{E}) precisely enumerates the possible trajectories of the
queue  giving rise to a given 2-TASEP configuration. The right
hand vector $|0\rangle$ represents an initial queue length of 0. At each service
time of the queue we have a matrix $D$ and at each non-service time a matrix
$E$. A vector $|n\rangle$ represents the length of the queue.  If the length of
the queue just before a service time is $n >0$ the action of $D$ on $|n\rangle$
is
\begin{equation}
D|n\rangle = |n\rangle + |n-1\rangle\;.
\label{Dac}
\end{equation}
The two terms represent the two possibilities at the service time:
the first  represents the service of a new  arrival at that time,
the second  represents   a service and no new arrival.
If $n=0$,
\begin{equation}
D|0\rangle = |0\rangle 
\label{Dac0}
\end{equation}
which  implies  that a new arrival has to be serviced at this time,
otherwise there would be an unused service which is forbidden.

Similarly, if the length of the queue is $n \geq 0$ 
just before a non-service time,
the action of $E$ on $|n\rangle$ is
\begin{equation}
E|n\rangle = |n\rangle + |n+1\rangle\;.
\label{Eac}
\end{equation}
The first term represents no new arrival  at that time,
the second term represents  a new  arrival  at that time.

The projector $|0 \rangle \langle  0 |$ at the left end of the string ensures
that only trajectories of the queue which finish at length 0 are counted
and the queue length is set to 0 for the start of the next string.\\

\noindent{\bf Remarks}
\begin{enumerate}
\item An alternative way to determine the queue length $n$ at a given time
  $t(i)=L-i$ is the following.  In Figure \ref{fig:2sq} each particle in line 1
  is associated to a particle in line 2 by a dashed black line. The length of
  the queue at $t(i)$ is given by the number of dashed black lines intersecting
  a vertical segment passing through $i^-$, i.e. just to the left of site~$i$;
  vertical dashed black lines do not affect the queue length.

\item The trajectory of the length of the queue in the queueing process 
with constraints described above  is  precisely a Motzkin path.
This allows one to represent matrix product calculations in terms of ensembles
of Motzkin paths see e.g. \cite{DS05,BCEPR06,BJJK04}.  
\end{enumerate}

\section{The $\eN$-Species TASEP}

 \begin{figure}[th]
  \begin{center}
\includegraphics[height=2.8cm]{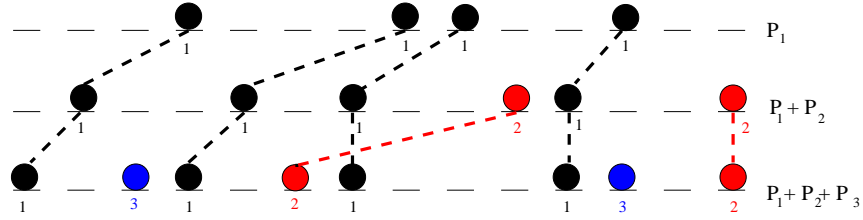}
   \caption{Graphical representation of Ferrari-Martin's algorithm.}
  \label{fig:algoFM}
\end{center}
\end{figure}

\subsection{Construction for $\eN$ Species of Particle}
In this section we review how the construction for the 2-species case
is extended to the $\eN$-species case \cite{FM07}.  For $\eN$ species of
particle we consider $\eN$-line configurations of particles.  The first
line comprises $P_1$ particles distributed randomly (with at most one
particle per site). The second line comprises $P_1 + P_2$ particles
distributed randomly and so on until the $\eN$th line which comprises
$P_1 + P_2 \cdots + P_\eN$ particles distributed randomly.  Initially,
in this $\eN$-line configuration, particles are not differentiated into
species. In the following we define the construction by which a
species label is attributed to each of the particles. Once this has
been done the $\eN$th line is identified with an $\eN$-TASEP
configuration.

We start with line 1 and associate each particle in line 1 to a particle in line
2 as in Section \ref{sec:angel}. This is done by beginning with a particle in
line 1 and associating it with the nearest particle, at the same site or to the
left, in line 2. We then take the next particle to the left in line 1 and
associate it in the same way to the first unassociated particle in line 2. This
process is continued until each of the particles in line 1 is associated with
one particle in line 2. These $P_1$ particles in line 2 are then labelled 1 and
the remaining $P_2$ unassociated particles in line 2 are labelled 2. The
resulting labels do not depend on which particle we began with in line 1, as
commented in Section \ref{sec:angel}.

We now proceed to associate the particles in line 2 with those in
line 3. First we use the same procedure as described above to
associate the particles labelled 1 in line 2 each to a particle in
line 3. The $P_1$ associated particles in line 3 are then
labelled 1. We then proceed to associate the particles labelled 2 in
line 2 to $P_2$ of the unassociated particles in line 3
(ignoring  the particles already labelled 1 in  line 3).
These particles in line 3 are then labelled 2 and the remaining $P_3$
unassociated particles in line 3 are labelled 3.

The procedure is then repeated up to line $\eN$ and results in $P_K$ of the
particles in line $\eN$ having label $K$ where $1\leq K \leq \eN$. The construction
is illustrated by an example in the three species case in Fig.
\ref{fig:algoFM}. Starting from the random distributions of particles in the $\eN$
lines, one obtains a configuration of the $\eN$ species TASEP with $P_K$ particles
of species $K$. The probability of the $\eN$-TASEP configuration so
  obtained is equal to its stationary probability under the $\eN$-TASEP
  dynamics. This was proven in \cite{FM07}; in  Section~\ref{sec:hat} we shall
  give an alternative proof, based on the matrix formulation, for the 3-TASEP.

Our aim is now to invert this construction. That is, for a given $\eN$-TASEP
configuration, we wish to compute the probability that it is generated by the
above procedure. This amounts to the combinatorial problem of counting the
number of ways particles may be distributed in the $\eN$-line configuration such
that the construction will lead to the desired $\eN$-TASEP configuration.

As has been discussed in Section \ref{sec:angel}, for the 2 species case Angel
\cite{Angel06} gave such a construction and this is equivalent to the matrix
product approach of Derrida {\it et al} \cite{DJLS93} (see section
\ref{sec:mp2}).  We showed this by using the queueing representation of Ferrari
and Martin.  In the following we first provide an algorithm, generalising the
pushing procedure of Section~\ref{sec:2} by which the probabilities can be
computed. We then construct explicit matrices which compute the
  $\eN$-TASEP weights by book-keeping the generalized pushing procedure.

\subsection{Ferrari and Martin's Construction and the Reverse Algorithm}
\label{sec:FM}
As described above, in the $\eN$-species procedure of Ferrari and Martin a
configuration of the $\eN$-TASEP is obtained from an $\eN$-line configuration, where
each of the $\eN$ lines consists of a single species TASEP configuration. The
procedure may be viewed in the following way: from line 1 (a configuration of
the single species TASEP) and line 2 one obtains a uniquely defined
configuration of the 2-TASEP; from that configuration of the 2-TASEP and line 3
one constructs a configuration of the 3-TASEP and so on, until one reaches a
configuration of the $\eN$-TASEP. Therefore, a given configuration of the
$\eN$-TASEP arises from a whole set of ($\eN -1$)-TASEP configurations that we shall
call its ancestors; each of these ($\eN -1$)-TASEP configurations arises in turn
from a whole set of ($\eN -2$)-TASEP ancestors etc...
  
 The stationary weight of the
initial $\eN$-TASEP configuration is then  given (but for an overall
normalization constant) by the sum of the weights of the ($\eN -1$)-TASEP
configurations that lead to it (its ancestors).  Applying this procedure
recursively, we observe that this stationary weight is given by the sum of the
weights of its ($\eN -2$)-TASEP ancestor configurations.  Finally, because the
single species TASEP has a uniform steady state, the weight of any $\eN$-TASEP
configuration is nothing but the total number of configurations of the single
species TASEP from which it derives. Therefore, to calculate the weight of a
given $\eN$-TASEP configuration we must determine the total number of 1-TASEP
configurations that are its ancestors.

  In the following we shall give a recursive algorithm to determine
  all the $(\eN-1)$-TASEP ancestors of a given $\eN$-TASEP configuration. Iterating
  this algorithm it is possible to obtain the total number of 1-TASEP ancestors
  of a given $\eN$-TASEP configuration; this number corresponds to the stationary
  weight of the $\eN$-TASEP configuration. 

In order to simplify our discussion we shall first present this algorithm for
the case of a 3-TASEP configuration, i.e. for a string of
particles of classes 1, 2, 3  
and holes (denoted by 0). We start from an initial  3-TASEP configuration.

  \begin{enumerate}
   
  \item  Freeze the positions of the  2s and the 3s.  Construct all
 possible configurations obtained by pushing  the 1s through the
 holes towards the right, until they hit a 2 or a 3 (i.e. a 1 can cross
 neither a  2 nor a 3).
\hfill\break
 $\rightarrow$  This procedure leads to  many  3-TASEP configurations with 
 various positions of the 1s. From now on the sites
 where the 1s are located will be passive.

  \item  Keep the positions of  the 3s frozen and start moving the 2s.
  For each configuration obtained above, construct all
 possible configurations obtained by pushing  the 2s through the
 holes towards the right, until they  hit  a 3. Note  that the  sites
 occupied by 1s  are spectators  and the 2s hop
 over  them as if they do not exist.
\hfill\break
  $\rightarrow$    This procedure leads to  many  3-TASEP configurations 
 in which the  positions of the 1s and the 2s are fixed.

\item Replace all the 3s by holes. We thus have obtained the
    complete set of 2-TASEP ancestors of the initial 3-TASEP configuration we
    started with.
   
\item The stationary weight of the initial 3-TASEP configuration
    (up to a global normalization constant) is the sum of the weights of its
    2-TASEP ancestors. 

  \end{enumerate}

\noindent Let us illustrate the algorithm with the  explicit example of the
string 2103

\begin{enumerate}
\item By pushing 1s to the right we obtain the strings 2103, 2013
\item By pushing 2s to the right (through the 1s), from 2103 we obtain
      2103 and 0123 and from 2013 we obtain 2013 and 0213
\item We replace 3s by 0s to obtain  the strings
      2100, 0120, 2010, 0210
    \item The weight of string 2103 in the 3-TASEP in terms of 2-TASEP weights
      is given by \[ W(2103)\, =\, W(2100)+ W(0120)+W(2010)+W(0210),\]
      which we may calculate, for example by using the matrix ansatz for the
      2-TASEP (\ref{DE}--\ref{AE}) and (\ref{0210}), as $W(2103) = 3 + 1+2 +3
      = 9$.
\end{enumerate}

It is easy to generalise to the $\eN$-species case and compute
the weight of a $\eN$-TASEP configuration in terms of the weights of 
$(\eN{-}1)$-TASEP configurations
  \begin{enumerate}
 
  \item   Freeze the positions of the  species $2,\ldots,K$.  Construct all
 possible configurations obtained by pushing  the 1s through the
 holes towards the right  (a 1 cannot  cross any species of particle).

  \item   Now in turn for $K=2,\ldots, N-1$ push species $K$
to the right keeping  the positions of  species $K+1,\ldots,N$ frozen and 
with species $1,\ldots, K-1$ spectators.
i.e.   For each configuration obtained from step 1, construct all
possible configurations obtained by pushing  the 2s through the
holes to the right, allowing the 2s to hop over 1s; then push 3s to the right
allowing 3s to hop over 2s and 1s, and so on until
species $K-1$ have been pushed to the right, 
hopping over all other species except $K$.

\item In all of the $\eN$-TASEP configurations generated in step (ii), replace all
  the $\eN$s by holes.  We thus have obtained a whole set of $(\eN{-}1)$-TASEP
  configurations: this is the complete set of ancestors of the initial $\eN$-TASEP
  configuration we started with.
   
\item  The sum of the stationary weight of all these $(\eN{-}1)$-TASEP ancestor
  configurations gives the stationary weight of the $\eN$-TASEP configuration. 

  \end{enumerate}

We saw that `the pushing procedure' of Angel is naturally implemented
by the $D, E$ and $A$ matrices. The algorithm given above is also
based on recursive pushing procedures  and as we shall
show in section \ref{sec:ma3} can be encoded by a matrix ansatz; in
this case the matrices for the $\eN$-TASEP are built by using the
matrices for the $(\eN{-}1)$-TASEP as elements.

\subsection{Queueing Interpretation of $\eN$-Species Construction}
\label{nsq}
Ferrari and Martin also proposed a queueing interpretation for the multiline
construction. The $\eN$ lines of the $\eN$-species construction correspond to
$\eN -1$ queues. The first line represents arrival times to the queue 1. The
second line represents service times for queue 1.  We continue using
  the convention that the queue time runs from right to left, so that the time
  $t(i)$ associated to site $i$ is given by $t(i)=L-i$.  As we have seen in
section~\ref{2sq}, unused service times of queue 1 become second-class
particles. Then when the particles of line 2 have been labelled either first or
second-class, they represent the arrivals for queue 2. The arrivals are
distinguished into first and second-class and the queue is a priority queue: at
the service times (given by the particles in line 3) the highest priority
waiting customer is always serviced first. That is, in queue 2 first-class
arrivals are served before second-class arrivals. The ouput of these service
times then become the arrival times for queue 3 with unused service times in
queue 2 providing third-class arrivals to queue 3. This construction is iterated
until the particles in line $\eN -1$, labelled $1,\ldots \eN-1$ provide the
arrivals for queue $\eN -1$ and the particles in line $\eN$ provide the service
times for queue $\eN -1$. When the particles in line $\eN$ are labelled $1
\ldots \eN$ they become the output of queue $\eN -1$, which corresponds to the
$\eN$-TASEP configuration.

 \section{The Matrix Product Formulation}

 In this section, we show how  the recursive
construction for the $\eN$-TASEP,  described in Section~\ref{sec:FM},
can be encoded within the  matrix product
formulation, described in Section~\ref{sec:mp2}.

 \subsection{Definition of the  Matrix Product  Ansatz
 and Simple Examples}

The matrix ansatz \cite{DEHP93} provides a solution to the stationary master
equation of the $\eN$-TASEP (made explicit later in (\ref{markov})), as
follows. First consider non-commuting matrices $X_0, X_1, \ldots, X_\eN$, where
$X_K$ is associated to particles of class $K$ (in particular $X_0$ is associated
to holes, $X_1$ is associated to first-class particles etc...). The ansatz
represents the stationary probability  $\prob({\cal C})$ of a
  $\eN$-TASEP configuration ${\cal C}=(\tau_1,\dots,\tau_L)$ 
 (where   $\tau_i$ is equal to $K$ if site $i$ 
 is occupied by a particle of class $K$)   as a statistical
  weight $W( {\cal C})$ divided by a normalization $Z$ 
\begin{equation}
 \prob({\cal C}) = \frac{1}{Z}W({\cal C})
 \label{ansatz}
\end{equation}
where the weight is given by the trace of the product of $L$ matrices, as
follows
\begin{equation}
W({\cal C})={\rm Tr}(X_{\tau_1}...X_{\tau_L})\;.
\end{equation}
Here $X_{\tau_i}$ is equal to $X_K$ if site $i$ is occupied by a particle of
class $K$ $(K=0,1,\ldots, N)$ in configuration ${\cal C}$.  The normalization
factor $Z$ (that depends on $L$ and on all the $P_K$'s where $P_K$ represents
the total number of particles of class $K$) ensures that $\sum_{\cal C} \prob({\cal C}) =
1$.  We emphasize that the matrix formulation depends on the number of species.
For example the matrices that represents first-class particles in the 2-TASEP
and the 3-TASEP are not the same.

If the system contains only first-class particles and holes, it is well known
that the stationary measure is uniform. 
Thus the particles and the holes may both be represented by one (a scalar) and
the matrix ansatz reduces here to a trivial form.

For the 2-TASEP holes, first-class and second-class particles are represented
respectively by the matrices $X_0=E$, $X_1=D$ and $X_2= A$ (in the notation of
\cite{DEHP93,DJLS93}) which satisfy (\ref{DE},\ref{DA},\ref{AE}).  It is
convenient to introduce matrices $\epsilon$ and $\delta$ defined by
\begin{equation}
      E = {\bf 1} + \epsilon,  \,\,\, 
   \,\,\,  D  = {\bf 1} + \delta  \, \, 
 \label{MPA1}
\end{equation}
where ${\bf 1}$ is the identity matrix. 
Then,  by (\ref{DE},\ref{DA},\ref{AE}), the matrices $\epsilon, A, \delta$
generate the following quadratic algebra:
 \begin{eqnarray}
   \delta\epsilon &=& {\bf 1}
   \nonumber \\
   \delta A  &=& 0 \nonumber \\
   A \epsilon  &=& 0  \;.
 \label{quadr}
 \end{eqnarray}
\subsection{Matrix Ansatz for 3-TASEP}\label{sec:ma3}
We now present  the matrix product formulation of  the stationary state
of the 3-TASEP.
\begin{eqnarray}
    X_1 &=&   {\bf 1}\otimes{\bf 1}\otimes D  + 
   \delta \otimes  \epsilon  \otimes A +   \delta \otimes {\bf 1}\otimes E
\label{MP31} \\
    X_2 &=&  A  \otimes {\bf 1}\otimes A
        +    A  \otimes  \delta \otimes E  
\label{MP32} \\
     X_3 &=&    A  \otimes A    \otimes E    
\label{MP33} \\
      X_0 &=&     {\bf 1}\otimes{\bf 1}\otimes E +
           {\bf 1}\otimes  \epsilon  \otimes A + 
   \epsilon  \otimes  {\bf 1}\otimes D \, .
\label{MP30} 
    \end{eqnarray}
    Note that the matrices are generally sums of tensor products of three
    semi-infinite matrices used in the matrix product representation of the
    stationary state of the 2-TASEP.  From the usual representation of the
    matrices $A, \delta$ and $\epsilon$ given in \cite{DEHP93,DJLS93} we obtain
    explicit expressions for the above matrices: they all have a block
    structure that is bidiagonal. Defining  matrices $F$, $G$, $H$
      and $\kK$  as

\begin{eqnarray}
F=               \left( \begin{array}{ccccc}    
                  D&0&0&0&\dots\\
                  0&D&0&0&\ddots\\
                  0&0&D&0&\ddots\\
                  0&0&0&D&\ddots\\              
                  \vdots &\ddots&\ddots&\ddots&
                        \end{array}  \right)
\qquad
G=
              \left( \begin{array}{ccccc}    
                  E&0&0&0&\dots\\
                  A&E&0&0&\ddots\\
                  0&A&E&0&\ddots\\
                  0&0&A&E&\ddots\\
                  \vdots &\ddots&\ddots&\ddots&
                        \end{array}  \right) \\[1ex]
H=               \left( \begin{array}{ccccc}    
                  A&E&0&0&\dots\\
                  0&A&E&0&\ddots\\
                  0&0&A&E&\ddots\\
                  0&0&0&A&\ddots\\              
                  \vdots &\ddots&\ddots&\ddots&
                        \end{array}  \right)
\qquad
\kK =
              \left( \begin{array}{ccccc}    
                  E&0&0&0&\dots\\
                  0&0&0&0&\ddots\\
                  0&0&0&0&\ddots\\
                  0&0&0&0&\ddots\\
                  \vdots &\ddots&\ddots&\ddots&
                        \end{array}  \right)
\end{eqnarray}
the matrix representation of (\ref{MP31}--\ref{MP30}) reads
\begin{eqnarray}
\hspace{-0.in} 
 X_1 =             \left( \begin{array}{ccccc}    
                  F&G&0&0&\dots\\
                  0&F&G&0&\ddots\\
                  0&0&F&G&\ddots\\
                  0&0&0&F&\ddots\\
                  \vdots &\ddots&\ddots&\ddots&
                        \end{array}  \right)
\qquad
  X_2 =               \left( \begin{array}{ccccc}    
                  H&0&0&0&\dots\\
                  0&0&0&0&\ddots\\
                  0&0&0&0&\ddots\\
                  0&0&0&0&\ddots\\
                  \vdots &\ddots&\ddots&\ddots&
                        \end{array}  \right)\\[1ex]
 X_3 =               \left( \begin{array}{ccccc}    
                  \kK &0&0&0&\dots\\
                  0&0&0&0&\ddots\\
                  0&0&0&0&\ddots\\
                  0&0&0&0&\ddots\\
                  \vdots &\ddots&\ddots&\ddots&
                        \end{array}  \right)\qquad
 X_0 =             \left( \begin{array}{ccccc}    
                  G&0&0&0&\dots\\
                  F&G&0&0&\ddots\\
                  0&F&G&0&\ddots\\
                  0&0&F&G&\ddots\\
                  \vdots &\ddots&\ddots&\ddots&
                        \end{array}  \right)
\end{eqnarray}
 All these matrices are triply infinite dimensional because 
 their  coefficients  are themselves infinite dimensional matrices 
 with elements  $D$,  $A$ and $E$ (which are also 
  infinite dimensional  matrices).

\subsection{Matrices as Priority Queue Matrices}

We now explain in the case $\eN =3$ how these matrices may be
  obtained  from the $\eN$-species queueing interpretation of the $\eN$-line
configuration discussed in section~\ref{nsq}.  In this case we have a 3-line
configuration that represents two queues in tandem.
Queue 1 has one type of customer which are considered to be first-class: Line 1
represents the arrival times of (first-class) customers in queue 1 and line 2
gives the service times of queue 1. The particles of line 2 are labelled 1 or 2
according to whether a service time is used or unused. Once labelled, the
particles of line 2 become the arrival times to queue 2. Queue 2 is a priority
queue containing first and second-class customers: any first-class customer is
served before the second-class customers  waiting in the
  queue. The particles of line 3 are the service times for queue 2. They are
labelled by which class of customer is served; if a service is unused it is
labelled 3.

We now consider the possible trajectories of the queue system.  To do this we
require 3 integer counters $l,m,n$: $l$ is the number of first-class customers
waiting in queue 2; $m$ is the number of second-class customers waiting in queue
2; $n$ is the number of first-class customers waiting in queue 1 (i.e. the
length of queue 1). The three counters indicate the state of the
  system at each queue time 
   $t(i) =L-i$ (which runs from right to left). 
  The counters  are  therefore  indexed by the times $t(i)$
  associated to sites $i$, but we omit this in our notation.

\noindent {\bf Remark:} The values of the counters can be obtained directly from
figure \ref{fig:algoFM} as follows:  for each site $i$, the counter
  $l$ with index $t(i)$ represents the number of black dashed lines crossing a
  vertical segment passing through $i^-$  between lines $2$ and $3$; $m$
represents the number of red dashed lines crossing the same segment and $n$
represents the number of black dashed lines crossing the segment between lines
$1$ and $2$; vertical dashed lines are not counted at all. The queue counters do
not register (a) second class particles served at their arrival time, (b) first
class particles served in both queues at their arrival time and (c) unused
services in the third line. However, a 3-TASEP configuration and the
trajectories of the three queues uniquely determine the 3-line configuration
generating it. This implies that it is enough to enumerate the set 
of queues trajectories compatible with the 3-TASEP configuration we are
computing the weight of.

 \hfill\break

The queue counters $l,m,n$   can be represented by a state vector $|l\, m\,
n\rangle$
\begin{equation}
|l\,  m\, n\rangle \equiv |l\rangle \otimes |m\rangle \otimes | n\rangle
\end{equation}
where  $|l\rangle =0$ for $l<0$.
We show now that the matrix product using $X_0,X_1,X_2,X_3$
defined in (\ref{MP31}--\ref{MP30}) precisely enumerate 
the possible trajectories of the state of the tandem queues  giving
  rise to a given configuration $(\tau_1,\dots,\tau_L)$.

We  list the possible updates of the counters $l,m,n$ at a given site
(or time), according to the line 3 label of that site, i.e.
the site variable $\tau_i$ in the $\eN$-TASEP configuration.
Then from each  possible update of the queue lengths
we deduce the necessary action of   the matrices $X_i$, $i=0,1,2,3$,
on the state vector $|l\, m\, n\rangle$.
Finally we can check
from  the definitions (\ref{MP31}--\ref{MP30}) of
the actions of $D$,$E$,$A$  (\ref{Dac},\ref{Dac0},\ref{Eac})
that $X_i | l\, m\, n\rangle$ produces the required update of the queue
counters.

\begin{description}
\item[$\tau_i=3$] In this case there is an unused service in queue 2 which
  implies $l=m=0$. In queue 1 there may or may have not been an arrival
  therefore $n\to n$ or $n\to n+1$. Thus, the action of $X_3$ must be
 \begin{equation}
 X_3 |l\,  m\, n\rangle =
\delta_{l,0} \delta_{m,0} 
 |0\rangle \otimes |0\rangle \otimes \left[ | n\rangle + |n+1\rangle\right]
= A \otimes A \otimes E |l\,  m\, n\rangle 
\end{equation}
which recovers the matrix expression for $X_3$, (\ref{MP33}).

\item[$\tau_i=2$] In this case a second-class service occurs in queue
  2 which implies that the number of first-class customers $l=0$ and
  there is no first-class arrival in queue 2.  If there  were a
  second-class arrival in queue 2 so that $m \to m$, it would imply
  $n=0$ as there would have to be an unusued service in queue 1. On the other
  hand, if there were no second-class arrival at queue 2 so that $m\to
  m-1$ then there might  or might  not be a first-class arrival at queue 1
  and $n\to n$ or $n\to n+1$.  Thus, the action of $X_2$ must be 
 \begin{eqnarray}
   X_2 |l\,  m\, n\rangle &=
   \delta_{l,0} |0\rangle \otimes\left[\, 
     \delta_{n,0} |m\rangle \otimes |0\rangle
     + 
     |m-1\rangle \otimes \left[|n\rangle + |n+1\rangle\right]\, \right]\\
   &= \left[A\otimes {\bf 1} \otimes A + A \otimes \delta \otimes E \right] |l\,  m\, n\rangle 
\end{eqnarray}
which recovers the matrix expression for $X_2$, (\ref{MP32}).

\item[$\tau_i=1$] In this case a first-class service occurs in queue 2. If there
  is also a first-class arrival at queue 2 then $l \to l$, $m\to m$ and $n\to
  n-$1 or $n$ since there is a first-class service and possibly a first-class
  arrival at queue 1. If there is instead a second-class arrival at queue 2 (a
  second-class service in queue 1) then there must be no first-class customers
  in queue 1 and so $l \to l-1$, $m\to m+1$ and $n=0$. Finally, if there is no
  arrival to queue 2 then there is no departure from queue 1 and there may or
  may not be an arrival at queue 1. Therefore $l \to l-1$, $m\to m$ and $n\to n$
  or $n+1$.

Thus, the action of $X_1$ must be
\begin{eqnarray*}
X_1 | l\, m\, n\rangle &=&
 |l\rangle \otimes |m\rangle \otimes \left[ | n\rangle +  | n-1\rangle\right]
+ |l-1 \rangle \otimes  |m+1 \rangle \otimes | 0\rangle \delta_{n,0} \\
&& +
|l-1\rangle \otimes |m\rangle \otimes \left[\, | n\rangle +  | n+1\rangle\,\right]\\
&=&
 |l\rangle \otimes |m\rangle \otimes D | n\rangle + 
\delta |l\rangle \otimes \epsilon |m\rangle \otimes A | n\rangle
+
\delta |l\rangle \otimes |m\rangle \otimes E | n\rangle \\
&=&
\left[ {\bf 1}\otimes{\bf 1}\otimes D  + 
   \delta \otimes  \epsilon  \otimes A +   \delta \otimes {\bf 1}\otimes E\right]
|l\, m\, n \rangle
\end{eqnarray*}
which recovers the matrix expression for $X_1$, (\ref{MP31}).

\item[$\tau_i=0$]
In this case there is no service at  queue 2.
If there is  first-class  arrival at queue 2 
$l \to l+1$, $m\to m$ and $n\to n-1$ or $n$ since
there is a first-class service and possibly a first-class arrival at queue 1.
If there is instead a second-class arrival at queue 2 (a second-class service
in queue 1) then there
must be no first-class customers in queue 1 and so
$l \to l$, $m\to m+1$ and $n=0$.
Finally, if there is no arrival at queue 2 then there is no departure from
queue 1  and there  may or may not be an
arrival at queue 1. Therefore $l \to l$, $m\to m$ and $n\to n$ or $n+1$.
Thus, the action of $X_1$ must be
\begin{eqnarray*}
X_0 | l\, m\, n\rangle &=&
 |l+1\rangle \otimes |m\rangle \otimes \left[\, | n\rangle +  | n-1\rangle\,\right]
+ |l \rangle \otimes  |m+1 \rangle \otimes | 0\rangle \delta_{n,0} \\
&& +
|l\rangle \otimes |m\rangle \otimes \left[\, | n\rangle +  | n+1\rangle\,\right]\\
&=&
 \epsilon|l\rangle \otimes |m\rangle \otimes D | n\rangle + 
|l\rangle \otimes \epsilon |m\rangle \otimes A | n\rangle
+
 |l\rangle \otimes |m\rangle \otimes E | n\rangle \\
&=&
\left[  \epsilon\otimes{\bf 1}\otimes D  + 
  {\bf 1} \otimes  \epsilon  \otimes A +     {\bf 1}\otimes {\bf 1}\otimes E\right]
|l\, m\, n \rangle
\end{eqnarray*}
which recovers the matrix expression for $X_0$, (\ref{MP30}).

\end{description}

\subsection{Algebraic Proof of the Matrix Product Ansatz}
\label{sec:hat}
The matrix product ansatz may be proved independently of the queueing
representation in an algebraic way.  We shall use the technique of ``hat
matrices'' to prove the ansatz (see e.g., \cite{BE07,sandow,nikolaus} for more
details).

We first recall the stationarity condition to be satisfied.  The dynamics of the
system can be encoded in a Markov matrix $\mM$ of size $\Omega \times \Omega$
where $\Omega$ is the total number of configurations of the system.  The
coefficient $\mM({\cal C},{\cal C}')$ of this matrix  represents the rate of
transition from a configuration ${\cal C}'$ to a different configuration ${\cal
  C}$; $-\mM({\cal C},{\cal C})$ is the total rate of exit from a given
configuration ${\cal C}$.  (Notice that this is the transpose of
  the usual generator matrix used in probability.) Thus the stationary
  probabilities 
must satisfy the stationary master equation
 \begin{equation}
 \sum_{{\cal C}'} \mM({\cal C},{\cal C}')\prob({\cal C'})  = 0 \,\,\,\, .
\label{markov}
 \end{equation}
 Due to the local structure of the rules (\ref{eq:rules1},\ref{eq:rules2}), $\mM$ 
 can be written as a sum of
 local matrices that represent the transitions
 that take place at  a bond
  $(i,i+1)$
\begin{equation}
 \mM = \sum_{i=1}^{L} \mM_{i,i+1} \,\,\,\, .
\end{equation}
$\mM_{i,i+1}$ are $(N+1)^2\times (N+1)^2$ matrices whose off diagonal elements
$\mM_{i,i+1}(\tau_i \tau_{i+1};\tau_i' \tau_{i+1}')$ give the transition rate from
configuration $\tau'_i \tau'_{i+1}$ to $\tau_i \tau_{i+1}$ at the bond $i,i+1$,
and whose diagonal element $\mM_{i,i+1}(\tau_i \tau_{i+1};\tau_i \tau_{i+1})$
gives minus the total transition rate out of configuration $\tau_i
\tau_{i+1}$. 
Since the only transitions involved in the $N$-TASEP are  exchanges at a bond, we have
  \begin{eqnarray}
    \label{p12}
    \mM_{i,i+1}(JK,KJ) = -\mM_{i,i+1}(KJ,KJ) = 1,\;
    &\hbox{if }K\ge 1\hbox{ and } J>K \hbox{ or }J=0\,, \nonumber\\
    \mM_{i,i+1}(K'J',KJ) = 0,&\hbox{otherwise, } \nonumber
  \end{eqnarray}
 where here $K,J$ are indices  that  take values from 0 to $\eN$. 
When the steady state probabilities are written 
in the matrix product form (\ref{ansatz})
the local  matrix $\mM_{i,i+1}$  acts only on the $i$th and
 the  $(i+1)$th matrices in the product.
The stationarity condition (\ref{markov}) then may be written
\begin{equation}
  \sum_{i=1}^L {\rm Tr}(X_{\tau_1}\ldots X_{\tau_{i-1}}
  Y_{\tau_i,\tau_{i+1}} X_{\tau_{i+2}} 
  \ldots X_{\tau_L})=0
\label{matstat}
\end{equation}
where
\begin{equation}
 Y_{\tau_i,\tau_{i+1}}  = \sum_{\tau_i',\tau_{i+1}'}
\mM_{i,i+1}(\tau_i \tau_{i+1};\tau_i' \tau_{i+1}')
X_{\tau_{i}'} X_{\tau_{i+1}'}\;.
\end{equation}
That is,
  \begin{equation}
    \begin{array}{ll}
      Y_{KJ} =  - X_KX_J&\hbox{for all}\quad K \ge 1 \hbox{ and }
      J > K \hbox{  or }  J =0\,, \\
      Y_{JK}  = X_KX_J&\hbox{for all}\quad K \ge 1 \hbox{ and }
      J > K \hbox{  or }  J =0\,,  \\
       Y_{JJ}=0 &\hbox{for all} \quad J\;.
    \end{array}
\label{p46}
  \end{equation}

  The key point to prove the validity of the matrix ansatz is to show that
  $Y_{\tau_i,\tau_{i+1}}$ is a divergence-like term, i.e.  there
    exist matrices ${\hat X_{\tau}}$ such that
\begin{equation}
\label{p47}
Y_{\tau_i,\tau_{i+1}}
 =   X_{\tau_i} {\hat X_{\tau_{i+1}}}
  - {{\hat X_{\tau_i}}} X_{\tau_{i+1}} \,   .
\label{cancellation}
\end{equation}
Summation over $i$ leads to a global cancellation in (\ref{matstat}), proving
thereby that the stationarity condition (\ref{markov}) is
satisfied. 
Combining (\ref{p46},\ref{p47}), we obtain the conditions: 
\def\b#1{\hat #1}
\begin{eqnarray}
  X_K  X_J &=&  \b{X}_K X_J -  X_K\b{X}_J
  \quad\quad \hbox{for all}\quad K \ge 1 \hbox{ and }
   J  >  K
 \hbox{  or }  J =0 
  \,,  \label{eq:hat1}\\
  X_K  X_J &=&  X_J\b{X}_K -  \b{X}_J X_K
  \quad\quad  \hbox{for all}\quad K \ge 1 \hbox{ and }
   J  >  K
 \hbox{  or }  J =0 \,,  \label{eq:hat2}\\
  0  
  &=&  X_J \b{X}_J-  \b{X}_J  X_J  
  \quad\quad\;\,\hbox{for all}\quad J.
    \label{eq:hat3}
\end{eqnarray}

 For $N=2$, it turns out to be rather easy to solve the above equations
 (see e.g.  \cite{BE07}):  indeed,   one finds that
(\ref{eq:hat1}--\ref{eq:hat3}) may be satisfied by choosing ${\hat X_{\tau}}$ to
be scalars so that they {\em commute} with ${ X_{\tau}}$. Then (\ref{eq:hat3})
is immediately satisfied and (\ref{eq:hat1},\ref{eq:hat2}) reduce to 3
conditions
\begin{eqnarray}
X_1 X_0 &=& {\hat X_1} X_0 - {\hat X_0} X_1\\ 
X_1 X_2 &=& {\hat X_2} X_1 - {\hat X_1} X_2\\
X_2 X_0 &=& {\hat X_2} X_0 - {\hat X_0} X_2\;. 
\end{eqnarray}
Choosing ${\hat X_1}=+1$, ${\hat X_0}=-1$,
${\hat X_2}=0$ and
$X_1=D$, $X_0=E$,
$X_2=A$ recovers
(\ref{DE}--\ref{AE}).
 
For $N=3$ it turns out  that choosing
${\hat  X_{\tau}}$ to be scalars does not allow (\ref{eq:hat1}--\ref{eq:hat2})
to be satisfied. Thus, the proof rests
upon finding the $4$ matrices ${\hat X}_K$ for $K=0,1,2, 3$.  
We now write  explicit forms for the hat matrices
that fulfil the above relations when $X_K$ are given by 
(\ref{MP31}--\ref{MP30}): 
\begin{eqnarray}
    \b{X}_1 &=&   ({\bf 1}-\delta) \otimes {\bf
1}\otimes{\bf 1} 
       \nonumber \\
    \b{X}_2 &=&  -A  \otimes \delta \otimes {\bf 1}       \nonumber \\
     \b{X}_3 &=&    -A  \otimes A    \otimes {\bf 1}    \nonumber \\
      \b{X}_0 &=& -X_0    +(\epsilon-{\bf 1})\otimes {\bf 1}\otimes{\bf 1} 
\;. \label{MPA3hat}        
\end{eqnarray}
It remains to verify that relations
(\ref{eq:hat1},\ref{eq:hat2},\ref{eq:hat3}) are satisfied. Here we
check a few relations involving $X_1$ and $\b{X}_1$.  For $J=1$ the
rhs of (\ref{eq:hat3}) becomes
\begin{eqnarray*}
 X_1 \b{X}_1 -  \b{X}_1  X_1
&=&   ({\bf 1}-\delta)\otimes{\bf 1}\otimes D  + 
   \delta ({\bf 1}-\delta) \otimes  \epsilon  \otimes A +   \delta({\bf
1}-\delta) \otimes {\bf 1}\otimes E\\
&-&
  ({\bf 1}-\delta)\otimes{\bf 1}\otimes D 
  -({\bf 1}-\delta) \delta \otimes  \epsilon  \otimes A -   ({\bf
1}-\delta)\delta \otimes {\bf 1}\otimes E =0\;.
\end{eqnarray*}
 Thus (\ref{eq:hat3}) is satisfied in the case $J=1$.

Using relations (\ref{quadr}) we find
\begin{eqnarray}
X_1X_2 =
 A \otimes {\bf 1} \otimes DA + A \otimes \delta \otimes DE
=  A \otimes {\bf 1} \otimes A + A \otimes \delta \otimes(D+E)\;,
\nonumber
\end{eqnarray}
and
\begin{eqnarray*}
  \b{X}_1 X_2 -  X_1 \b{X}_2
&=&  A \otimes {\bf 1} \otimes A + A \otimes \delta \otimes E
  -(-A \otimes \delta \otimes D) \\
&=& X_1 X_2\\
  X_2\b{X}_1  -   \b{X}_2 X_1
&=&   A({\bf 1}-\delta)\otimes{\bf 1}\otimes A  + 
    A({\bf 1}-\delta)\otimes \delta \otimes E  +\\
&& -(- A \otimes \delta \otimes D    -A\delta \otimes \delta \epsilon\otimes A
-A \delta \otimes \delta \otimes E)\;,\\
&=& X_1 X_2
\end{eqnarray*}
thus  (\ref{eq:hat1},\ref{eq:hat2}) are satisfied for the case $K=1$, $J=2$.
Similarly, all relations
 (\ref{eq:hat1},\ref{eq:hat2},\ref{eq:hat3}) may be verified.

 \section{Hierarchical Matrix Ansatz for the Multispecies ASEP}
\label{sec:n>3}

In this section, we generalize the previous construction to the multispecies
totally asymmetric exclusion process 
on the ring ${\mathbb Z}_L$
with $\eN$ classes of particles for any
$\eN >1$.  We show that a matrix ansatz for a system containing $\eN$ classes of
particles (plus holes) can be constructed recursively knowing a matrix ansatz
for a system with $\eN -1$ classes of particles (plus holes).  We shall simply
present the results here and give some examples; we  leave the algebraic proof
and further generalizations to a future publication.

  The matrices   $X_K^{(\eN)}$
  at level  $\eN$  are obtained  by making tensor products
  of  the   $X_M^{(\eN-1)}$ defined at  level $\eN -1$ with 
 some   matrices  $a_{KM}$  constructed from 
  $\epsilon, A, \delta$  and ${\bf 1}$.
 The  matrix ansatz is  given by
\begin{eqnarray}
        X_0^{(\eN)} &=&    \sum_{M=0}^{\eN-1}  a_{0M}^{(\eN)} \otimes X_M^{(\eN-1)}
 \label{def:X0}    \\ 
      X_K^{(\eN)}  &=&  a_{K0}^{(\eN)}  \otimes   X_0^{(\eN-1)}
 +  \sum_{M=K}^{\eN-1}  a_{KM}^{(\eN)} \otimes   X_M^{(\eN-1)} \, \, 
\,\,\,   \hbox{  for } 1 \le K \le \eN\;.
 \label{def:XK}
\end{eqnarray}
We emphasize that in this section our notation for the matrix $X_K^{(\eN)}$ has
two indices: the lower index $K$ denotes the class of the particle represented
by the matrix, whereas the upper index $\eN$ gives the total number of classes
considered in  the system.

The fundamental building blocks to construct the $a_{KM}^{(\eN)}$
matrices are the matrices $\epsilon, A, \delta$ and ${\bf 1}$
(identity).  The $a_{KM}^{(\eN)}$ are then given by
\begin{eqnarray}
   a_{00}^{(\eN)} &=&  {\bf 1}^{ \otimes (\eN-1) }    \label{def:a00}  \\
   a_{0M}^{(\eN)} &=&   {\bf 1}^{ \otimes (M-1) }  \otimes  \epsilon
    \otimes {\bf 1}^{ \otimes (\eN-M-1) }
  \,\,\,   \hbox{  for } 1 \le M \le \eN-1\;.
  \label{def:a0M}
\end{eqnarray}
 For $K \ge 1$ we have 
\begin{eqnarray}
  a_{K0}^{(\eN)}  &=&  A^{ \otimes (K-1) } \otimes  \delta 
   \otimes {\bf 1}^{ \otimes (\eN-K-1) }     \,\,\,   \hbox{  for }
     1 \le K \le \eN-1
  \label{def:aK0}   \\ 
 a_{KK}^{(\eN)}  &=&   A^{ \otimes (K-1) } \otimes 
 {\bf 1}^{ \otimes (\eN-K) }  \label{def:aKK}   \\ 
  a_{KM}^{(\eN)}  &=&  A^{ \otimes (K-1) } \otimes  \delta 
   \otimes {\bf 1}^{ \otimes (M-K-1) }  \otimes  \epsilon
    \otimes {\bf 1}^{ \otimes (\eN-M-1) }  
 \,\,\,    \hbox{  for } 1 \le     K  <  M \le \eN-1 {\hskip 0.75cm}
  \label{def:aKM} \\
 a_{\eN 0}^{(\eN)}  &=&   A^{ \otimes (\eN-1) }\;.   \label{def:aN0}    
\end{eqnarray}
 It is understood in the formulae above  that any  matrix
 raised to a tensor-power equal to zero is equal to the scalar 1
 which can  be removed from the tensor product.

 Note from (\ref{def:X0},\ref{def:XK}) that 
  the matrices   $X_K^{(\eN)}$  at level  $\eN$ are composed of tensor products
of ${{\eN} \choose 2}$ fundamental matrices $\epsilon, A, \delta$ or ${\bf 1}$. 

 \subsection{Some  Examples}

   Using   the  hierarchical matrix ansatz given above, we study
 explicitly  the cases
 $\eN  \le 3$. 

   For $\eN =0$, the system does not contain any particles but only holes.
 There is only one configuration which has probability 1. Thus, 
 we   {\it define}  $X_0^{(0)} =1$.

  For $\eN =1$, we obtain from equations~(\ref{def:X0}) and~(\ref{def:XK}),
 using the fact that  $X_0^{(0)} =1$
 \begin{equation}
   X_0^{(1)} =  a_{00}^{(1)} \,\,\, \hbox{ and } \,\,\,
    X_1^{(1)} =  a_{10}^{(1)} \;.
 \end{equation}
 But from equations~(\ref{def:a00}) and~(\ref{def:aN0}) we 
 find that $ a_{00}^{(1)} =  a_{10}^{(1)}  = 1$ and we recover 
 the fact  that for a system with only  one class of particles
 the stationary measure is uniform and therefore the matrix ansatz
 is trivial.

 For $\eN =2$, we find from equations~(\ref{def:a00}), (\ref{def:a0M}),
 (\ref{def:aKK}) and (\ref{def:aN0}), that $a_{00}^{(2)} = {\bf 1}$,
 $a_{01}^{(2)} = \epsilon$, $ a_{10}^{(2)} = \delta $, $ a_{11}^{(2)} = {\bf
   1}$, and $ a_{20}^{(2)} = A$. Then, from the recursion
 relations~(\ref{def:X0}) and~(\ref{def:XK}), and using the fact that $X_0^{(1)}
 = X_1^{(1)} = 1$, we deduce that
\begin{eqnarray}
  &X_0^{(2)} =  a_{00}^{(2)} + a_{01}^{(2)} =  {\bf 1} +  \epsilon =E,  
  \,\,\, 
  \,\,\,   X_1^{(2)} =  a_{10}^{(2)} +  a_{11}^{(2)} = {\bf 1} +  \delta = D, 
  \,\,\,\, \label{ma201}\\
  \,\,\,& \hbox{ and }    X_2^{(2)} = a_{20}^{(2)}  =  A \, .\label{ma22}
 \end{eqnarray}
 We retrieve the  fundamental matrix ansatz~(\ref{MPA1}).

For $\eN =3$,  using the recursion relations
~(\ref{def:X0}) and~(\ref{def:XK}), the definitions  
(\ref{def:a00}--\ref{def:aN0}) and the results (\ref{ma201}),(\ref{ma22}), we obtain
\begin{eqnarray}
X_0^{(3)} &=&    a_{00}^{(3)}\otimes X_0^{(2)}+
a_{01}^{(3)}\otimes X_1^{(2)}+a_{02}^{(3)}\otimes X_2^{(2)}\nonumber \\
&=&
 {\bf 1}\otimes{\bf 1}\otimes E +   \epsilon  \otimes  {\bf 1}\otimes D +
  {\bf 1}\otimes  \epsilon  \otimes A \\
X_1^{(3)} &=&    a_{10}^{(3)}\otimes X_0^{(2)}+
a_{11}^{(3)}\otimes X_1^{(2)}+a_{12}^{(3)}\otimes X_2^{(2)}\nonumber \\
  &=&  \delta \otimes {\bf 1}\otimes E +  {\bf 1}\otimes{\bf 1}\otimes D  + 
  \delta \otimes  \epsilon  \otimes A   \\
X_2^{(3)} &=&    a_{20}^{(3)}\otimes X_0^{(2)}+
a_{22}^{(3)}\otimes X_2^{(2)} =   A  \otimes  \delta \otimes E  + A  \otimes {\bf 1}\otimes A  \\
X_3^{(3)} &=&    a_{30}^{(3)}\otimes X_0^{(2)}
          = A  \otimes A    \otimes E\, , 
    \end{eqnarray}
and we  retrieve the  expressions (\ref{MP31}--\ref{MP30}).

\section{Discussion}
In this work we have considered the multispecies totally asymmetric exclusion
process on the ring ${\mathbb Z}_L$ (although our results are generalisable to
${\mathbb Z}$).  We have shown how the stationary measure may be written in a
matrix product formulation, thus providing an algebraic proof of the stationary
measure which we presented for the three species case $\eN =3$. For arbitrary
$\eN$ we have shown how the matrix product formulation may be constructed in a
hierarchical fashion, although we leave the algebraic proof of the stationary
measure to a future publication.

 Ferrari and Martin have constructed the
  stationary state of the $\eN$-TASEP as the output of $N$ queues in series with
  $N$ priority-classes of customers, for all $\eN$. The construction takes a
  $\eN$-line binary configuration sampled at random and produces a $\eN$-TASEP
  configuration whose resulting law is invariant for the $\eN$-TASEP.  We have
  shown that the matrix ansatz for $\eN =2$ of Derrida et al \cite{DJLS93} gives
  a mechanism to count the number of $2$-line configurations producing a given
  $2$-TASEP configuration; the matrices may be thought of as acting on the space
  of queue counters. For $\eN =2$ there is just a single queue with one type of
  customer and the queue counter is simply the length of the queue. We have also
  extended the matrix ansatz for $\eN>2$.  In this case there are multiple
priority queues and there are several queue counters representing the number of
each class of customer in each queue. This results in the queue matrices acting
on tensor product spaces and accordingly the `matrices' of the matrix
formulation become higher rank tensors.  This relation between the matrices and
queueing processes also provides us with a natural representation of the space
on which the matrices act. Until now, it was believed that the matrices act on a
purely formal `auxiliary' space which did not have any physical interpretation.

The algebraic proof of the stationary measure for $\eN >3$ relies on the
existence of `hat' matrices \cite{sandow,nikolaus,BE07} described in
Section~\ref{sec:hat}.  This is in contrast to the $\eN =2$ case where the hat
matrices were simply scalars and the relations obeyed by the matrices become a
quadratic algebra, as in the $\eN=1$ open-boundaries case of \cite{DEHP93}.  For
$\eN >2$ the relations between the matrices have a more complicated algebraic
structure and it would be of interest to explore this further.

One advantage of the matrix product formulation of the stationary
measure is that it provides a framework within which the calculation of
quantities of physical interest, such as correlation functions, can be
carried out. So far we have not attempted such calculations but it would be 
important to do so.

Finally, we mention that the multispecies TASEP may be generalised in several
ways by introducing rates which differ from one or additional processes.  For
example, for $\eN =1$ allowing particles to carry out forward exchanges with
holes with rate $p$ and backward exchanges with rate $q$ generates the partially
asymmetric exclusion process for which a matrix product formulation of the
steady state  on the open boundary system  has been fully worked
out \cite{Sasamoto,BECE}.  In the case $\eN =2$ there are partially asymmetric
generalisations which admit matrix product stationary states \cite{DJLS93,BE07}.
So far, in work in progress, we have found a partially asymmetric generalization
of the matrix product ansatz presented in Section~\ref{sec:n>3}.  It would be of
interest to understand  how the queueing interpretation of the steady state should
be modified.

\section*{Acknowledgements}
KM thanks Nikolaus Rajewsky for inspiring discussions and many memorable moments
devoted to the $N$-TASEP model.

We thank the Isaac Newton Institute, Cambridge for hospitality during the
programme {\em Principles of Dynamics of Nonequilibrium Systems} where this work
was begun.  MRE thanks the CNRS for a Visiting Professorship and the Laboratoire
de Physique Th\'eorique et Mod\`eles Statistiques, Universit\'e Paris-Sud
for  hospitality.

\end{document}